\documentclass[12pt]{amsart}
\usepackage{amsmath,amssymb,amsfonts,amsthm}
\usepackage{graphicx}
\usepackage{enumerate}
\usepackage{tikz}
\usepackage[mathscr]{eucal}
\usepackage{float}
\numberwithin{equation}{section}
\newtheorem{theorem}{Theorem}[section]

\newtheorem{corollary}[theorem]{Corollary}

\numberwithin{equation}{section}

 \makeatletter
      \def\@setcopyright{}
      \def\serieslogo@{}
      \makeatother

\begin{document}
\bibliographystyle{plain}

\title{Short Paper No. 2: Relations between integer sequences A003622 $\&$ A035336 and two criteria in respect of maximum degrees within the family of finite linear Jaco graphs}
\author{Timoth$\acute{e}$e Raso$^{\dag}$ and Johan Kok$^{\ddag}$} 
\address{$^{\dag}$Independent Researcher and Engineer, Lyon, France.}
\email{
timothee.raso@gmail.com}
\address{$^{\ddag}$Independent Mathematics Researcher, City of Tshwane, South Africa.}
\email{jacotype@gmail.com}
\date{}
\begin{abstract}
This paper establishes two relationships between integer sequences and the maximum degrees of finite linear Jaco graphs $J_n(x)$. First, it characterizes the positive maximum degrees that occur exactly once within the family, showing that they are generated by the integer sequence A003622. Second, it determines the orders of the corresponding graphs and establishes their relationship with the integer sequence A035336. These results resolve conjectures previously stated in \cite{kok}. The proofs were developed with the assistance of AI-based mathematical reasoning technology and subsequently reviewed by the authors.
\end{abstract}
\maketitle
\section{Preliminary remarks}
Several conjectures concerning integer sequences and graph parameters were proposed in \cite{kok}. The first author subsequently investigated these conjectures using AI-based mathematical reasoning technology under development. The resulting arguments were reviewed by the second author, and both authors collaborated on the preparation of this manuscript. Further details concerning the use of AI-assisted reasoning are provided in Section \ref{sec1}.

Some definitions and background material from \cite{kok1, raso1} are recalled to make the present paper self-contained. The results established here concern the multiplicity of maximum degrees and the orders of the corresponding finite linear Jaco graphs, extending the maximum degree formula established in \cite{raso1}. The present manuscript is intended for submission to a peer-reviewed mathematical journal.
\section{Introduction}
We assume familiarity with standard graph theoretic terminology and notation; see \cite{bondy, harary, west}. For a graph $G$, its maximum degree is denoted by $\Delta(G) = max\{deg_G(v): v \in V(G)\}$. The cardinality of a set $X$ is denoted by $|X|$. Unless otherwise stated, the graphs considered in this paper are the underlying graphs of finite directed linear Jaco graphs, referred to here simply as finite linear Jaco graphs; see \cite{kok1, kok2}.

The study is motivated by two conjectured relationships between integer sequences and graph parameters.
Building on the known formula for the maximum degree of $J_n(x))$, we investigate which positive maximum degrees occur exactly once and determine the orders of the corresponding graphs. The resulting proofs also provide an example of the use of AI-assisted mathematical reasoning in the investigation of mathematical conjectures.
\section{On linear Jaco graphs}
Since our interest lies with the underlying graph we recall the construction method as explained in \cite{raso1}.  

\textit{Construction method to obtain $J_n(x)$, $n \geq 4$:} 

(a) Let $J_1(x) = P_1$, $J_2(x) = P_2$ and $J_3(x) = P_3$.

(b) Let the empty graph (edgeless graph) of order $m \geq 1$ be denoted by $\mathfrak{N}_m$. Let $X = \{v_i: i = 4,5,6,\dots,n\}$

and let $G_2 = P_3 \cup \mathfrak{N}_{n-3}$ where $V(\mathfrak{N}_{n-3}) = X$.

(c) By normal consecutive step-count for $i = 3,4,5,\dots, n-1$ do as follows: 

To obtain $G_i$ add the edges $v_iv_{i+1}$, $v_iv_{i+2}$,$\dots, v_iv_{i+t}$ with $t$ a maximum such that, $deg_{G_i}(v_i) \leq i$.

(d) After completion of step-count $i = n-1$ label the resultant graph $J_n(x)$.
Note that $J_n(x)$ is a simple connected undirected graph.

Recall from \cite{kok, raso1} that for a given $J_n(x)$ the open neighborhood of $v_i$ is given by:
$N(v_i) = N^-(v_i) \cup N^+(v_i)$ where $N^-(v_i) = \{v_j:j < i$ and the edge $v_jv_i$ exists$\}$ and $N^+(v_i) = \{v_j:j > i$ and the edge $v_iv_j$ exists$\}$.
Also let $t_1(v_i) = |N^-(v_i)|$ and $t_2(v_i) = |N^+(v_i)|$ and let $t^\star_2(v_i) = max\{t_2(v_i)$ in $J_n(x)$ evaluated at $v_i$: $n \geq i\}$. The aforesaid yields $t^\star_2(v_i) = i - t_1(v_i)$. We recall an important result from \cite{raso1}.

\begin{theorem}\label{thm1}\cite{raso1}
For every $n \geq 1$ and every $1 \leq i \leq n$ and every $v_i \in V(J_n(x))$,
\begin{center}
$deg(v_i) = min\{i, n-\lfloor \frac{i+1}{\varphi}\rfloor\}.$
\end{center}
\end{theorem}
\subsection{Non-zero, non-repetitive maximum degree of certain $J_n(x)$, $n = 1,2,3,\dots$:} 
Recall from \cite{raso1} that:
\begin{theorem}(Maximum degree of each $J_n(x)$, $n = 1,2,3,\dots$)\label{thm2}\cite{raso1}
For every integer $n \geq 1$ we have,
\begin{center}
$\Delta(J_n(x)) = \lfloor \frac{n+2}{\varphi}\rfloor -1 = \lfloor \frac{2(n+2)}{1 + \sqrt{5}}\rfloor -1$
\end{center}
where the golden ratio is given by $\varphi = \frac{1+\sqrt{5}}{2}$.
\end{theorem}
For $n = 1,2,3,\dots$ it is easy to generate a sequence of maximum degrees of $J_n(x)$, namely
\begin{center}
$0,1,2,2,3,3,4,5,5,6,7,7,8,8,9,10,10,11,11,12,13,13,14,15,15,\dots$
\end{center}
Now a strictly increasing sequence of non-zero, non-repetitive maximum degrees of certain finite linear Jaco graphs can be extracted i.e.
\begin{center}
$1,4,6,9,12,14,\dots$
\end{center}
The first question to be settled is, whether there exists a closed formula to generate the strictly increasing sequence of the non-zero, non-repetitive maximum degrees within the family of finite linear Jaco graphs?

Furthermore, a strictly increasing sequence representing the order $n \geq 2$ of the corresponding finite linear Jaco graphs can be extracted i.e.
\begin{center}
$2,7,10,15,20,23,\dots$
\end{center}
The second question is, whether there exists a closed formula to generate the strictly increasing sequence of the order $n \geq 2$ of the corresponding finite linear Jaco graphs?
\section{AI-assisted mathematical reasoning:} The known identities which will be used are,
\begin{center}
$\varphi^2 = \varphi +1$, $\frac{1}{\varphi} = \varphi-1$, $\frac{1}{\varphi^2} = 2-\varphi$, $\frac{1}{\varphi} + \frac{1}{\varphi^2} = 1.$
\end{center}
Both, Beatty and Wythoff sequences have been studied extensively. The reader is referred to \cite{ginosar, kimberling} for a worthy overview.
\subsection{Multiplicity of the maximum degrees of $J_n(x)$, $n \in \mathbb{N}$:}
\begin{theorem}\label{thm3}
For $t = 1,2,3,\dots$ the values given by sequence A003622 i.e. the Wythoff compound sequence AA: $a(t) = \lfloor t \times (\frac{1+\sqrt{5}}{2})^2\rfloor -1$ are the non-zero, non-repetitive ('unique') maximum degrees of certain (or corresponding) finite linear Jaco graphs.
\end{theorem}
\begin{proof}
Let the integer $d \geq 1$ and let $k = d+1$. From Theorem \ref{thm2} it follows that,
\begin{center}
$\Delta(J_n(x)) = d \Leftrightarrow \lfloor \frac{n+2}{\varphi}\rfloor = k$.
\end{center}
Because $k\varphi$ and $(k+1)\varphi$ are irrational it follows that $k\varphi < n+2 < (k+1)\varphi$. So for $d \geq 1$ we have $k \geq 2$. Hence, every integer within the interval $k\varphi < n+2 < (k+1)\varphi$ corresponds to $n \geq 2$. Therefore, the multiplicity of $d$ is
\begin{center}
$M(k) = \lfloor (k+1)\varphi\rfloor - \lfloor k\varphi\rfloor$.
\end{center}
Since $1 < \varphi < 2$, $M(k) \in \{1, 2\}$. Using $\varphi = 1 + \frac{1}{\varphi}$ we have 
\begin{center}
$M(k) = 1 + \underbrace{\lfloor \frac{k+1}{\varphi}\rfloor - \lfloor \frac{k}{\varphi}\rfloor}_{=A}$.
\end{center}
Now $A$ is exactly one when an integer $t$ lies strictly between $\frac{k}{\varphi}$ and $\frac{k+1}{\varphi}$ or put differently, when:
\begin{center}
$k < t\varphi < k+1$.
\end{center}
By the irrationality of $t\varphi$ it is equivalent to $k = \lfloor t\varphi\rfloor$. Therefore,
\begin{center}
\begin{equation}
M(k) = 2 \Leftrightarrow k \in \{\lfloor t\varphi\rfloor: t \geq 1\}.
\end{equation}
\end{center}

It follows from Beatty's theorem, namely if $\alpha, \beta > 1$ are irrational and $\frac{1}{\alpha} + \frac{1}{\beta} = 1$ then the sequences $\lfloor t\alpha\rfloor$ and $\lfloor t\beta\rfloor$, $t \geq 1$, partition the positive integers. Taking $\alpha = \varphi$ and $\beta = \varphi^2$ gives
\begin{center}
\begin{equation}
\mathbb{Z}_{>0} = \{\lfloor t\varphi\rfloor:t \geq 1\}\cup \{\lfloor t\varphi^2\rfloor: t \geq 1\}.
\end{equation}
\end{center}
Combining equations $(4.1)$ and $(4.2)$ gives $M(k) = 1 \leftrightarrow k = \lfloor t\varphi^2\rfloor$ for some $t \geq 1$. Since $k = d+1$ it follows exactly that,
\begin{center}
$|\{n \geq 2: \Delta(J_n(x)) = d\}| = 1 \Leftrightarrow d = \lfloor t\varphi^2\rfloor-1$.
\end{center}
or
\begin{center}
$|\{n \geq 2: \Delta(J_n(x)) = d\}| = 1 \Leftrightarrow d = \lfloor t(\frac{1 +\sqrt{5}}{2})^2\rfloor-1$, $t= 1,2,3,\dots$.
\end{center}
Note that $d = \lfloor t(\frac{1 +\sqrt{5}}{2})^2\rfloor-1$, $t= 1,2,3,\dots$ is exactly the Wythoff compound sequence AA: $a(t) = \lfloor t \times (\frac{1+\sqrt{5}}{2})^2\rfloor -1$, $t= 1,2,3,\dots$ found in \cite{sloane}. This settles the result.
\end{proof}
\subsection{Order of $J_n(x)$, $n \in \mathbb{N}$ which yields non-zero, non-repetitive maximum degrees:}
Note that the proof of Theorem \ref{thm3} gives the complete multiplicity dichotomy. If $k = d+1$ then
\begin{equation*} 
\mid\{n \geq 2: \Delta(J_n(x)) = d\}| = 
\begin{cases}
$1$ & \mbox{if $k = \lfloor t\varphi^2\rfloor$ for some $t \geq 1$},\\
$2$ & \mbox{if $k = \lfloor t\varphi\rfloor$ for some $t \geq 1$}.
\end{cases}
\end{equation*} 

Thus every non-zero positive maximum degree occurs either once or twice. The interval description $k\varphi < n+2 < (k+1)\varphi$ contains the information needed to determine the corresponding order $n$ whenever the multiplicity is one. The aforesaid forms the basis of this subsection. We identified an incorrect lower bound for $t$ in Conjecture 2.4 of \cite{kok}. The corrected bound is $t = 1,2,3,\dots$
\begin{theorem}\label{thm4}
Let $n \geq 2$. The value $\Delta(J_n(x))$ occurs exactly once in the sequence
\begin{center}
$\Delta(J_1(x)), \Delta(J_2(x)), \Delta(J_3(x)), \Delta(J_4(x)), \Delta(J_5(x)), \dots$
\end{center}
if and only if there is a unique integer $t \geq 1$ such that $n = 2\lfloor t\varphi\rfloor + t - 1$.
\end{theorem}
\begin{proof}
By Theorem \ref{thm3} a positive maximum degree $d$ occurs exactly once if and only if $d = \lfloor t\varphi^2\rfloor - 1$ for some integer $t \geq 1$. The parameter $t$ is unique because $\lfloor t\varphi^2\rfloor$ is strictly increasing. Let $k = d+1 = \lfloor t\varphi^2\rfloor$ and $a = \lfloor t\varphi\rfloor$. Since $\varphi^2 = \varphi + 1$ we have $k = a + t$. By Theorem \ref{thm2} the condition $\Delta(J_n(x)) = d$ is equivalent to $k\varphi < n+2 < (k+1)\varphi$.

By Theorem \ref{thm3} the interval $k\varphi < n+2 < (k+1)\varphi$ contains exactly one integer, namely $m = \lfloor k\varphi\rfloor + 1$. Write $t\varphi = a + x$, where $0 < x < 1$. Then $k\varphi = (a + t)\varphi = 2a + t + x(2 - \varphi)$.

Since $0 < x(2 - \varphi) < 1$ we have $\lfloor k\varphi\rfloor = 2a + t$. Consequently,
\begin{center}
$n + 2 = m = 2a + t + 1$
\end{center}
which gives,
\begin{center}
$n = 2\lfloor t\varphi\rfloor + t - 1$.
\end{center}
Conversely, for every integer $t \geq 1$ define
\begin{center}
$k = \lfloor t\varphi^2\rfloor$ and $n = 2\lfloor t\varphi\rfloor + t - 1$.
\end{center}
The  preceding calculation shows that $n + 2$ is the unique integer in the interval $(k\varphi, (k + 1)\varphi)$. Therefore, $\Delta(J_n(x)) = k - 1$ and Theorem \ref{thm3} guarantees that this maximum degree occurs exactly once. This settles the result.
\end{proof}
\begin{corollary}\label{cor1}
For every integer $t \geq 1$, the $t$-th term of the integer sequence A035336,
\begin{center}
$n(t) = 2\lfloor \frac{t}{\varphi}\rfloor + t - 1$
\end{center}
is the order of the unique finite linear Jaco graph whose maximum degree is equal to the $t$-th term of the integer sequence A003622,
\begin{center}
$d(t) = \lfloor \frac{t}{\varphi^2}\rfloor -1$.
\end{center}
Thus, the two sequences correspond term by term, with A003622 giving the positive maximum degrees of multiplicity one and A035336 giving the orders of their corresponding graphs.\footnote{See \cite{kok} for the corresponding $p$-graphical parts of Corollary \ref{cor1}.}
\end{corollary}
\section{AI disclosure}\label{sec1}   
AI-assisted mathematical reasoning technology was used substantially in the mathematical exploration, development of candidate arguments, intermediate verification, and preparation of the source manuscript. The technology was used within an author-directed problem-solving and verification methodology. All resulting mathematical claims and proofs were subsequently reviewed and validated before human drafting of a final manuscript. The human authors take full responsibility for all mathematical statements, proofs, computations, citations, and conclusions presented herein.
\section*{In Memory}
The concept of linear Jaco graphs followed from a dream during the night 10/11 January 2013 which was the first dream Kokkie (second author) could recall about his daddy, after his daddy passed away in the peaceful morning hours of 24 May 2012, shortly before the awakening of Bob Dylan, celebrating Dylan's 71 st birthday. The dream involved the following problem:
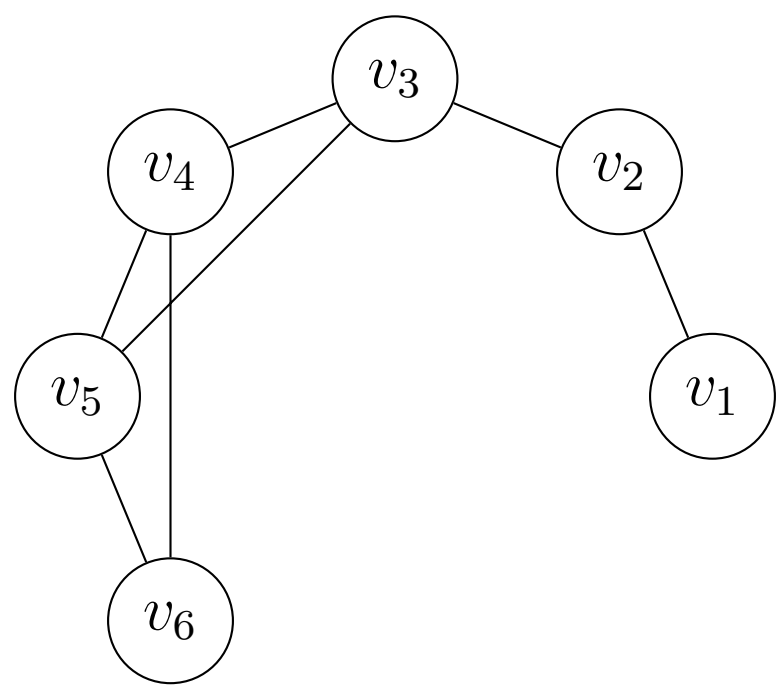
\begin{figure}[h]
\centering
\begin{tikzpicture}[scale=0.7]
\tikzstyle{every node}=[draw,shape=circle];
\node (v0) at (0:3) [] {$v_1$};
\node (v1) at (45:3) [] {$v_2$};
\node (v2) at (2*45:3) [] {$v_3$};
\node (v3) at (3*45:3) [] {$v_4$};
\node (v4) at (4*45:3) [] {$v_5$};
\node (v5) at (5*45:3) [] {$v_6$};
\draw (v0) -> (v1)
(v1) -- (v2)
(v2) -- (v3)
(v2) -- (v4)
(v3) -- (v4)
(v3) -- (v5)
(v4) -- (v5)
;
\end{tikzpicture}
\caption{The Dream}\label{fig1}
\end{figure}

\end{document}